\documentclass[twoside]{article}
\usepackage{graphicx}
\usepackage{caption}
\usepackage[utf8]{inputenc}
\usepackage{subcaption}
\usepackage{amsmath}  %% ADDED by FM for \binom
\title{\bf  %A note on the 
\hskip 0.2truecm Wright functions of the second kind \\ and Whittaker functions}
\author{\sc F.Mainardi$^1$,  R.B. Paris$^2$,  and  A. Consiglio$^3$ 
\\
{\em $^1$Dipartimento di Fisica e Astronomia, Universit\`{a} di Bologna, \& INFN,}\\
{\em Via Irnerio 46, I-40126 Bologna, Italy}\\
{\em  E-mail: francesco.mainardi@bo.infn.it}
\\
{\em $^2$Division of Computing and Mathematics, University of Abertay, Dundee DD1 1HG, UK}\\
{\em  E-mail: r.paris@abertay.ac.uk}
\\
{\em $^3$Institut f\"{u}r Theoretische Physik und Astrophysik and W\"{u}rzburg-Dresden Cluster of}\\
{\em Excellence ct.qmat, Universit\"{a}t W\"{u}rzburg, 97074 W\"{u}rzburg, Germany}\\
{\em  E-mail: armando.consiglio@physik.uni-wuerzburg.de}
\\
}
\begin{document}
\def\f#1#2{\mbox{${\textstyle \frac{#1}{#2}}$}}
\def\dfrac#1#2{\displaystyle{\frac{#1}{#2}}}
\def\boldal{\mbox{\boldmath $\alpha$}}
\def\ds{\displaystyle}
\def\cen{\centerline}
\def\e{\hbox{e}}
\def\Ai{\hbox{Ai}} 
\def\d{\partial}
%%%%%%%%%%%%%%%% DEFINITIONS of NUMBER SETS  %%%%%%%%%%
\def\RR{\vbox {\hbox to 8.9pt {I\hskip-2.1pt R\hfil}}}
\def\NN{{\rm I\hskip-2pt N}}
\def\CC{{\rm C\hskip-4.8pt \vrule height 6pt width 12000sp\hskip 5pt}}
%%%%%%%
%{\newcommand{\Sgoth}{S\;\!\!\!\!\!/}
\newcommand{\bee}{\begin{equation}}
\newcommand{\ee}{\end{equation}}
\newcommand{\la}{\lambda}
\newcommand{\ka}{\kappa}
\newcommand{\al}{\alpha}
\newcommand{\fr}{\frac{1}{2}}
\newcommand{\fs}{\f{1}{2}}
\newcommand{\g}{\Gamma}
\newcommand{\br}{\biggr}
\newcommand{\bl}{\biggl}
\newcommand{\ra}{\rightarrow}
\newcommand{\gl}{\raisebox{-.8ex}{\mbox{$\stackrel{\textstyle >}{<}$}}}
\newcommand{\gtwid}{\raisebox{-.8ex}{\mbox{$\stackrel{\textstyle >}{\sim}$}}}
\newcommand{\ltwid}{\raisebox{-.8ex}{\mbox{$\stackrel{\textstyle <}{\sim}$}}}
\renewcommand{\topfraction}{0.9}
\renewcommand{\bottomfraction}{0.9}
\renewcommand{\textfraction}{0.05}
\newcommand{\mcol}{\multicolumn}
\date{}
\maketitle
\pagestyle{myheadings}
\markboth{\hfill \sc F. Mainardi, R. B.\ Paris and  A. Consiglio   \hfill}
{\hfill \sc  Wright functions of the second kind and Whittaker functions\hfill}
%\vskip -0.25truecm
%\cen{\bf arXiv version 05.11.21 MPC-WRIGHT-WHITTAKER-FINAL.tex} 
% \cen{\bf  with 3 Sections and Appendices: A:Paris' details, B:Inversion of LT}
%%
%\cen{\bf Vol 24, No 1, pp. 54--72 (2021) DOI: 10.1515/fca-2021-0003} 
\vskip 0.3truecm

\begin{abstract}

In the framework of higher transcendental functions the Wright functions of the second kind
have increased their relevance resulting from their applications
in probability theory and, in particular, in fractional diffusion processes.
Here, these functions are compared with the well-known Whittaker functions in some special
cases of fractional order.
In addition, we point out two erroneous representations in the literature.

%%%%%%%%%% Enter suitable ky words and phrases:
\noindent
{\bf Keywords}: {Fractional calculus, Wright functions, Whittaker functions,
 Hypergeometric functions,  Laplace transform.}

\noindent
{\bf Mathematics Subject Classification}:{26A33, 30B10, 30E15, 33C20, 34E05,  41A60}

\noindent
{\bf Paper published in Fract. Calc. Appl. Anal. Vol. 25 (2022), pp. 858-875.}
\\
{\bf DOI: 10.1007/s13540-022-00042-2}

\end{abstract} %%%%%%%%%%%%%

%%%%%%%% begin papers' body %%%%%%%%%%%%%%%%%%%%%%%%%%%%%

%%%%%%%%%%%%%%%%%%%%%%%%%%% Section 1 %%%%%%%%%%%%%%%%%%%
\section{Introduction}

\setcounter{section}{1}\setcounter{equation}{0}

The  Wright function under consideration (also known as a generalised Bessel function) is defined by
\bee
\label{e10}
\phi(\lambda, \mu, z) :=  W_{\lambda,\mu}(z)=\sum_{n=0}^\infty\frac{z^n}{n! \g(\lambda n+\mu)}, \quad z  \in \CC,
\ee
where $\lambda$ is supposed real and $\mu$ is, in general, an arbitrary complex parameter. The series
converges for all finite $z$ provided $\lambda>-1$ and, when $\lambda=1$, it reduces to the modified Bessel function $z^{(1-\mu)/2}I_{\mu-1}(2\sqrt{z})$.
We point out that the $\phi$ notation is originally due to Wright while the $W$ notation was introduced by Mainardi
{in the 1990's, see \cite{Mainardi_WASCOM1993,Mainardi_AML1996}}.

For the Wright function corresponding  to a negative $\lambda$ it is convenient to denote by $\nu$ the positive parameter
 $ -\lambda =\nu$, with $0<\nu<1$.
The function with negative $\lambda$ has been termed a Wright function of the second kind by Mainardi \cite{FM1}, 
with the function with $\lambda>0$ being referred to as a Wright function of the first kind.

In order to avoid confusion with the identical notation for the Whittaker function, we shall denote in the following
 this latter function by $\mathcal{W}$.
Hereafter we recall the definition of the Whittaker functions, which are confluent hypergeometric functions as  found
in the NIST Handbook \cite[(13.2.42), (13.14.5)]{DLMF}.
Indeed, the Whittaker function ${\cal W}_{\kappa,\mu}(z)$ can be expressed in terms of the confluent hypergeometric function ${}_1F_1(z)$ by
\begin{equation} \label{e18a}
\begin{array}{ll}
{\cal W}_{\kappa,\mu}(z)=e^{-z/2} z^{\frac{1}{2}+\mu}&\left\{\dfrac{\Gamma(-2\mu)}{\Gamma(\frac{1}{2}-\mu-\kappa)}\,{}_1F_1(\fs+\mu-\kappa;1+2\mu;z)\right. \\ \\ % 
&\left. +\, \dfrac{z^{-2\mu}\Gamma(2\mu)}
{\Gamma(\fs+\mu-\kappa)}\,{}_1F_1(\fs-\mu-\kappa;1-2\mu;z)\right\},
\end{array}
\end{equation}
where
\[
{}_1F_1(a;b,z)=\sum_{k=0}^\infty \frac{(a)_k z^k}{(b)_k k!} \quad (|z|<\infty).
\]
For more details the reader is referred,
for example, to the NIST Handbook \cite[p.~334]{DLMF}.

 \medskip %%%%
 
The plan of our paper is as follows.
In Section \ref{sec:2}
 we give the definition of the auxiliary functions $F_\nu(z)$ and $M_\nu(z)$ as special cases of the Wright function $W_{-\nu,\mu}(z)$ with $0<\nu<1$ and $0\leq\mu\leq1$. Various plots of these functions for real argument are shown to illustrate their behaviour. The relation of these function to the Whittaker function is indicated.
 In Section \ref{sec:3}
 we present a list of special evaluations of $W_{-\nu,\mu}(\pm x)$ for certain rational values of $\nu$ in the range $0<\nu<1$ and $0\leq\mu\leq1$ expressed in terms of Whittaker, Airy and Bessel functions.
  Section \ref{sec:4} deals with a Laplace transform pair arising in time fractional diffusion processes 
 and its relation to the so-called ``four sisters''. 
 The proof of a typical result stated in Section \ref{sec:3} %%
  is given in Appendix A. The proof of the Laplace transform pair discussed in Section 
  \ref{sec:4}  is given in Appendix B.
 Section \ref{sec:5} is devoted to concluding remarks.

%%%%%%%%%%%%%%%%%%% Section 2 %%%%%%%%%%%%%%%%%%%%%%%%%%%

\section{The Wright function of the second kind versus the Whittaker function} \label{sec:2}

\setcounter{section}{2} \setcounter{equation}{0}

The Wright functions arise in probability theory related to the analysis of some Levy-stable distributions\footnote{We recall on this respect the representations of the stable distributions obtained by 
Schneider \cite{Schneider LNP1986}
in terms of Fox $H$-functions and in some cases in terms of Whittaker functions that would be
compared with ours.}  
and, more specifically, in processes governed by time-fractional diffusion and diffusion-wave
equations. Indeed,  partial differential equations of non-integer order in time
 \begin{equation} 
 \dfrac{\d^\beta u}{\d t^\beta} = \dfrac{\d^2 u}{\d x^2},
 \quad u= u(x,t), \  \; \beta = 2\nu \in (0,2)\,,
\label{Fractional PDE}
\end{equation} 
 were outlined by Mainardi in the early 1990's; see, for example, \cite{MT}.
 For more details, see the 2010 book by Mainardi \cite{FM1}, the recent survey by Consiglio and Mainardi
\cite{Consiglio-Mainardi  NLFO21},  and references therein.

 In the above context, the following {\it auxiliary functions} were introduced:
\bee\label{e11a}
F_\nu(z)=W_{-\nu, 0}(-z)=\sum_{n=1}^\infty \frac{(-z)^n}{n! \g(-n\nu)},
\quad 0<\nu<1,
\ee
\bee\label{e11b}
M_\nu(z)=W_{-\nu,1-\nu}(-z)=\sum_{n=0}^\infty \frac{(-z)^n}{n! \g(-n\nu+1-\nu)},
\quad 0<\nu<1.
\ee
These functions are interrelated by the following relation:
\bee\label{e12}
F_\nu(z)=\nu z M_\nu(z).
\ee
For the asymptotic expressions of the Mainardi auxilary functions we refer the reader to
\cite{MT,Mainardi-Tomirotti GEO97}, \cite{PCM FCAA21}.
 For further information {\bf about} the general Wright functions we refer to the papers by
 Luchko \cite{LL,Luchko HFCA19} and by Paris \cite{P14,P17,PHbk} and references therein.

For particular rational values of the parameter $\nu$ the Wright functions of the second kind 
are expected to be represented in terms of known special functions of  the hypergeometric class.
For instance, referring to the $M$-Wright functions with positive variable $x$
the following particular representations are nowadays well known in terms of some simpler functions:
\bee\label{e13}
\begin{array}{ll}
&M_0(x)=\e^{-x},
\quad M_{1/3}(x)=3^{2/3} \mbox{Ai}(x/3^{1/3}),
\quad M_{1/2}(x)=\frac{1}{\sqrt{\pi}}\,\e^{-x^2/4},\\ \\
&M_{2/3}(x) = 3^{-2/3}
\left[3^{1/3}\,x\, \Ai\left(x^2/3^{4/3}\right) -
3\Ai^\prime\left ( x^2/3^{4/3}\right)\right ]
\, \e^{-2x^3/27},
\end{array}
\ee
where Ai is the Airy function and Ai$'$ its derivative. As $\nu \to1^-$ the function 
$M_\nu(x)$ tends to the Dirac delta generalized function $\delta(x-1)$.

\begin{figure}[h!]
	\centering
	\includegraphics[width=.35\textwidth,angle=+0.0]{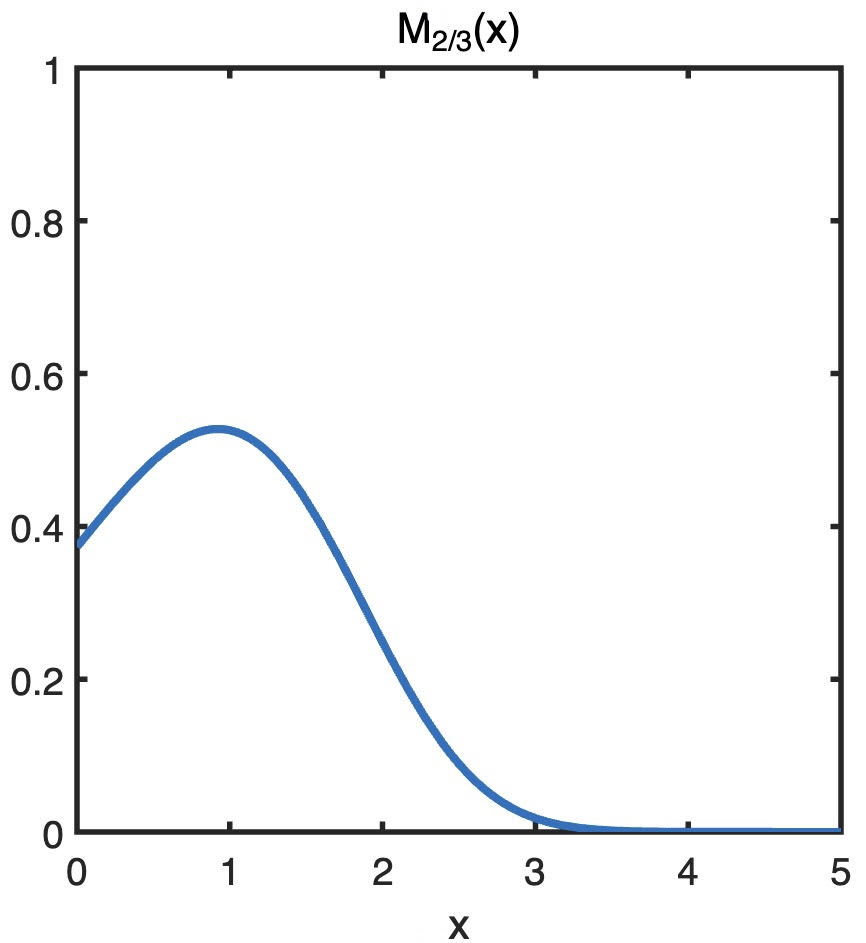}
\vskip -0.1truecm
		\caption{The function $M_{2/3}(x) = W_{-2/3, 1/3}(-x)$ for
 $0\le x\le 5$}
 \label{fig:1} 
	\end{figure}

 Plots of $M_\nu(|x|)$ for real $x$ and varying $\nu$ are presented in \cite[Appendix F]{FM1} and \cite{Mainardi-Consiglio SI20} in order to illustrate the transition between the special values $\nu=0, 1/2, 1$, 
 and the physical transition between diffusion and wave propagation.
As an example we include in Fig.~\ref{fig:1} %% 1 %% 
the plot of the function $M_{2/3}(x) = W_{-2/3, 1/3}(-x)$ for  $0\le x\le 5$.

%%%%%%%%%%%%%% Figure 1 %%%%%%%%%%%%%%%%%%%%%%%%%%%%%
  %%%%%%%%%%%%%%%%%%%%%%

For details on the $M$-Wright function in probability theory, see for example the papers by Cahoy \cite{Cahoy CS2012} and by Pagnini \cite{Pagnini FCAA2013}.
 The case $\mu=0$ in (\ref{e10}) also finds application in probability theory and is discussed extensively by Paris and Vinogradov in \cite{PV} where it is referred to as a `reduced' Wright function.
In both representations we write
\bee\label{e12a}
\phi(\lambda,0;z)=W_{\lambda,0}(z)\,.
\ee
The simple representations (\ref{e13}) have motivated us to explore the possibility that for rational values of
$\nu= 1/4, 1/3, 1/2, 2/3, 3/4$ the Wright functions of the second kind can be represented
in terms of other special functions (of hypergeometric type) including the Whittaker functions.
In particular, we concentrate our attention on the representations with  $\nu = 2/3$,
where we are able to correct two  erroneous results  existing in the literature.
We note that another method for representation of the Wright function in terms of the
hypergeometric functions for the rational values of its parameters was introduced
in the paper by Gorenflo et al. 
(see Section 2.2 of \cite{Gorenflo-et-al_FCAA1999}).
%\cite{Gorenflo-et-al_FCAA1999}.
This method is based on the representation of the Wright function as a particular case of the Fox $H$-function and on using the Gauss-Legendre formula for the gamma function.
% (see Section 2.2 of \cite{Gorenflo-et-al_FCAA1999}).
\newpage

\begin{figure}[h!]
	\centering
	\includegraphics[width=.75\textwidth,angle=+0.0]{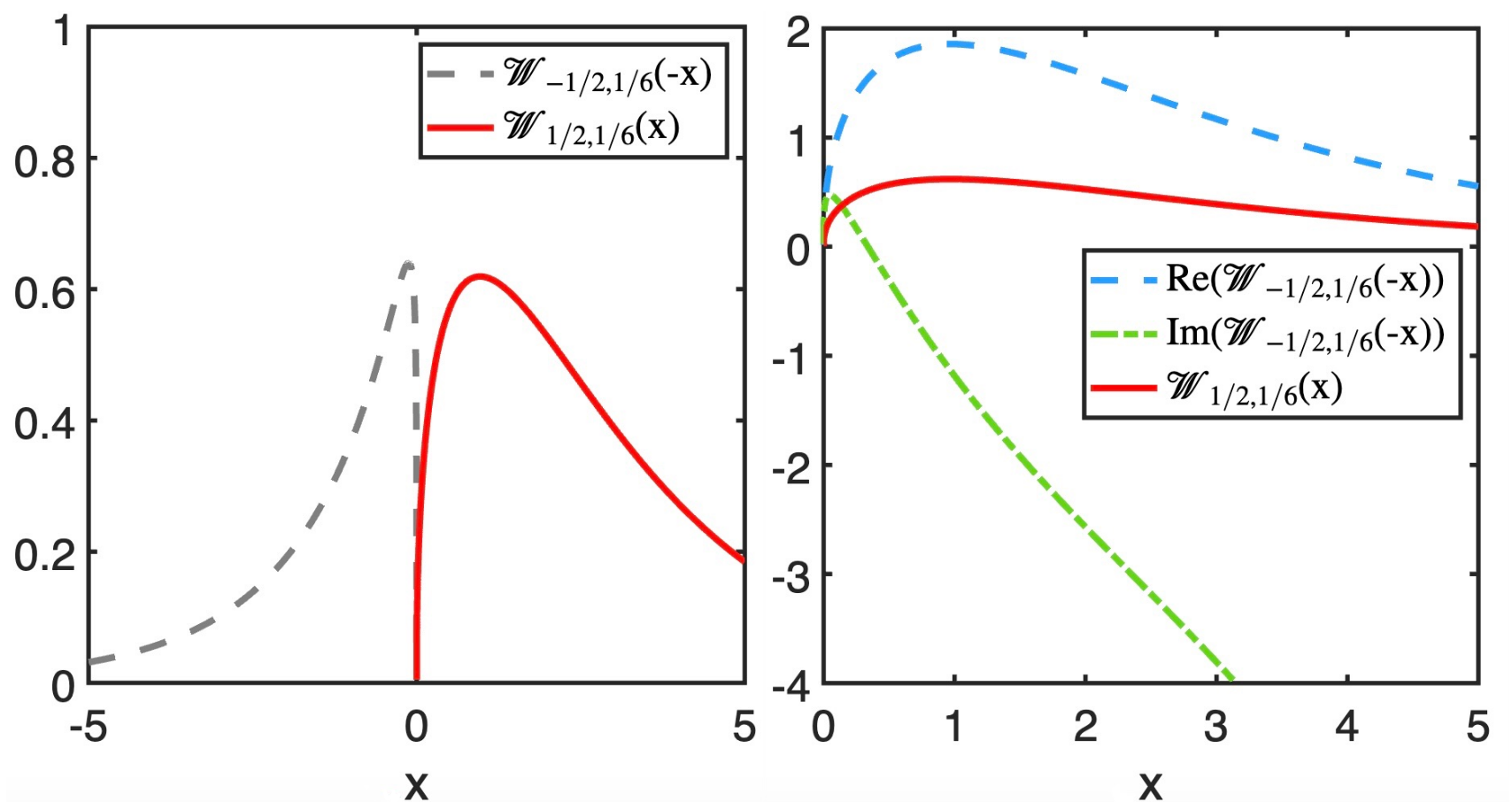}
	%\includegraphics[width=.75\textwidth,angle=+0.0]{CONSIGLIO_W_TWO.png}
%%	\vskip -0.3truecm %%%	
	\caption{The two solutions
 ${\mathcal W}_{1/2, 1/6}(x), \, {\mathcal W}_{-1/2, 1/6}(-x)$
 of the Whittaker equation with $0\le x\le 5$.}
 \label{fig:2} %%%
	\end{figure}

We recall that the Whittaker functions are so named after the fundamental 1903 paper
by Whittaker \cite{Whittaker 1903}.
They are particular confluent hypergeometric functions that are the solutions of the
following differential equation
\begin{equation} 
\label{Whittaker equation}
\frac{d^2}{d x^2} {\mathcal W}_{\mu,\nu}(x) + \left(- \frac{1}{4} +
  \frac{\mu}{x} + \frac{1/4-\nu^2} {x^2} \right) {\mathcal W}_{\mu,\nu}(x) = 0.
\end{equation}
In \cite{Whittaker 1903} Whittaker noted (using our notation)
that the differential equation  is unchanged if
$ \nu$  is replaced by $- \nu$  and
if $\mu$  is replaced by $- \mu$, provided $x$ is replaced by
$- x$  at the same time. Hence the four functions
${\mathcal W}_{\mu,\nu}(x) $,
${\mathcal W}_{\mu,-\nu}(x) $,
${\mathcal W}_{-\mu,\nu}(-x) $,
${\mathcal W}_{-\mu,-\nu}(-x) $
are solutions of the differential equation (\ref{Whittaker equation}).
Then we have two linearly independent solutions of the Whittaker equation (\ref{Whittaker equation})
$ 
{\mathcal W}_{\mu,\nu}(x), \;
{\mathcal W}_{-\mu,\nu}(-x),\; x\ge 0
$;
 see \cite[p.~6]{Erdelyi-Swanson BOOK57} and the NIST Handbook \cite[p.~335]{DLMF}.
 We note that the Whittaker functions exhibit a branch cut on the negative real axis, so that
they assume complex values on this semi-axis.
 This fact is clarified  in  Fig.~\ref{fig:2} %% 2 %% 
 concerning the above Whittaker functions for
$0\le x\le 5$ in   the special cases $\mu=\pm 1/2$ and $\nu=1/6$ of most interest in the following.

%%%%%%%%%%%%%% Figure 2 %%%%%%%%%%%%%%%%%%%%%
 %%%%%%%%%%%%%%%%%

	In \cite{Stankovic 70}, Stankovi\'{c}, in addition to wrongly  reporting  the
Whittaker differential equation (\ref{Whittaker equation}), which was presumably a misprint,
he obtained the following representation of the Wright function (also reported in the treatise on Mittag-Leffler functions by Gorenflo {et al.} \cite{GKMR Book20} in Eq. (7.2.12), p. 214)
$$ 
W_{-\frac{2}{3}, 0} (- x^{-\frac{2}{3}}) =  - \frac{1}{2 \sqrt{3\pi}}
\exp\left(- \frac{2}{27 x^2}\right)  {\mathcal W}_{- \frac{1}{2}, \frac{1}{6}}\left(- \frac{4}{27 x}\right).
$$  
We note that Stankovi\'{c}'s  representation appears to be wrong because the corresponding Whittaker function 
is expected to be complex valued.

In order to derive the correct result we take advantage of the
 Whittaker function representations for the reduced Wright functions
$W_{-2/3,0}(\pm x)$ given by Paris and Vinogradov in {\cite[Appendix C]{PV}
for $ x\ge 0$ and checked in the plots in Fig.~\ref{fig:C6-C7}, %%% 3 %% 
\bee
\phi(-2/3,0,x)=W_{-2/3,0}(x) = - \dfrac{1}{2\sqrt{3\pi}}
\e^{2x^{3}/27}\,
{\mathcal{W}}_{-1/2, 1/6}\left(\dfrac{4x^3}{27}\right),
\label{C6}
\ee
\vskip -0.3truecm %% 
 \bee
\phi(-2/3,0,-x)=W_{-2/3,0}(-x) =  \sqrt{\dfrac{3}{\pi}}
\e^{-2x^{3}/27}\,
{\mathcal{W}}_{1/2,1/6}\left(\dfrac{4x^3}{27}\right),
\label{C7}
\ee

%%%%%%%%%%%%%% Figure 3 %%%%%%%%%%%%%%%%%%%
\begin{figure}[h!]
\centering
\includegraphics[width=.40\textwidth,angle=+0.0]{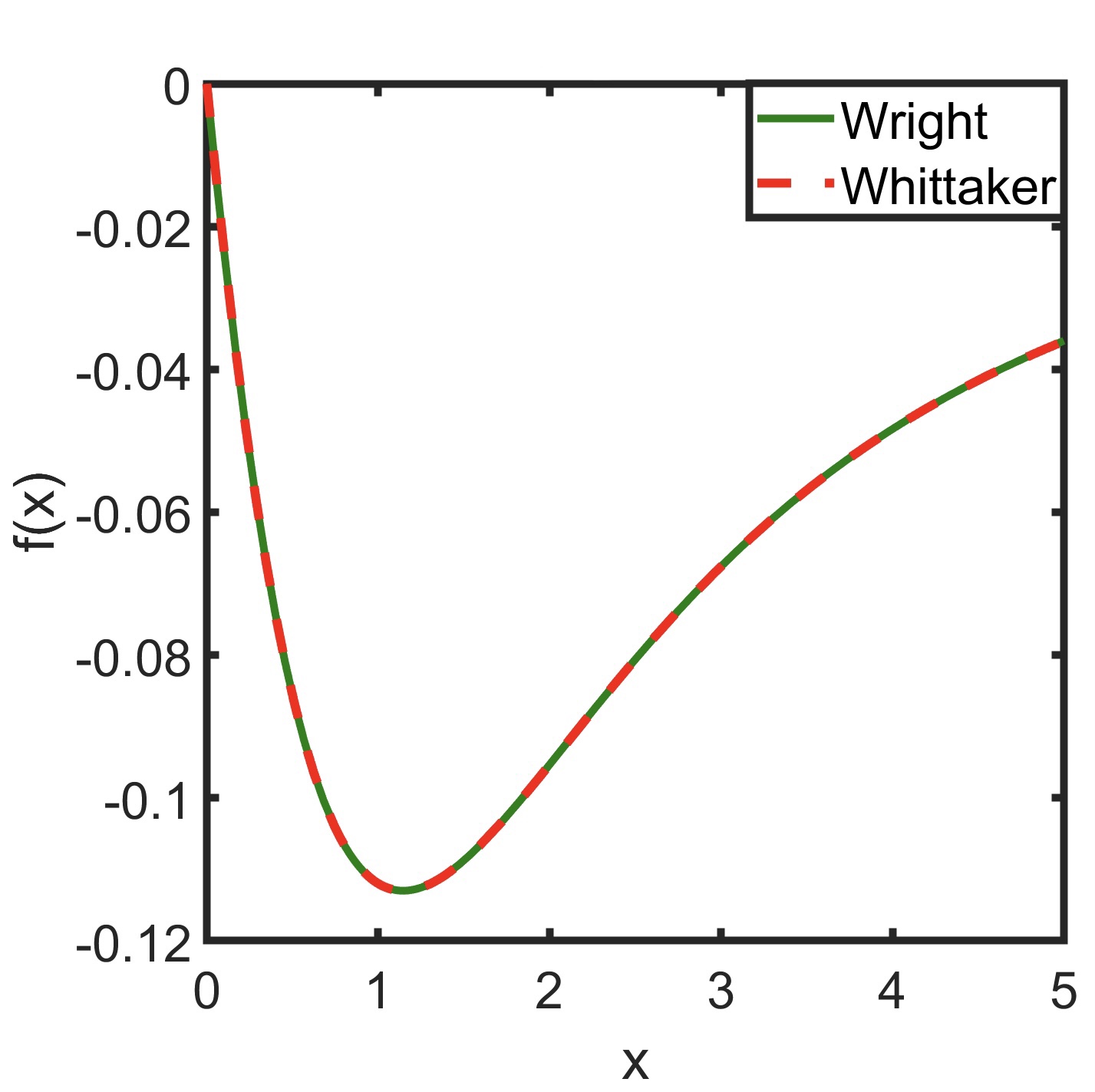}
\includegraphics[width=.38\textwidth,angle=+0.0]{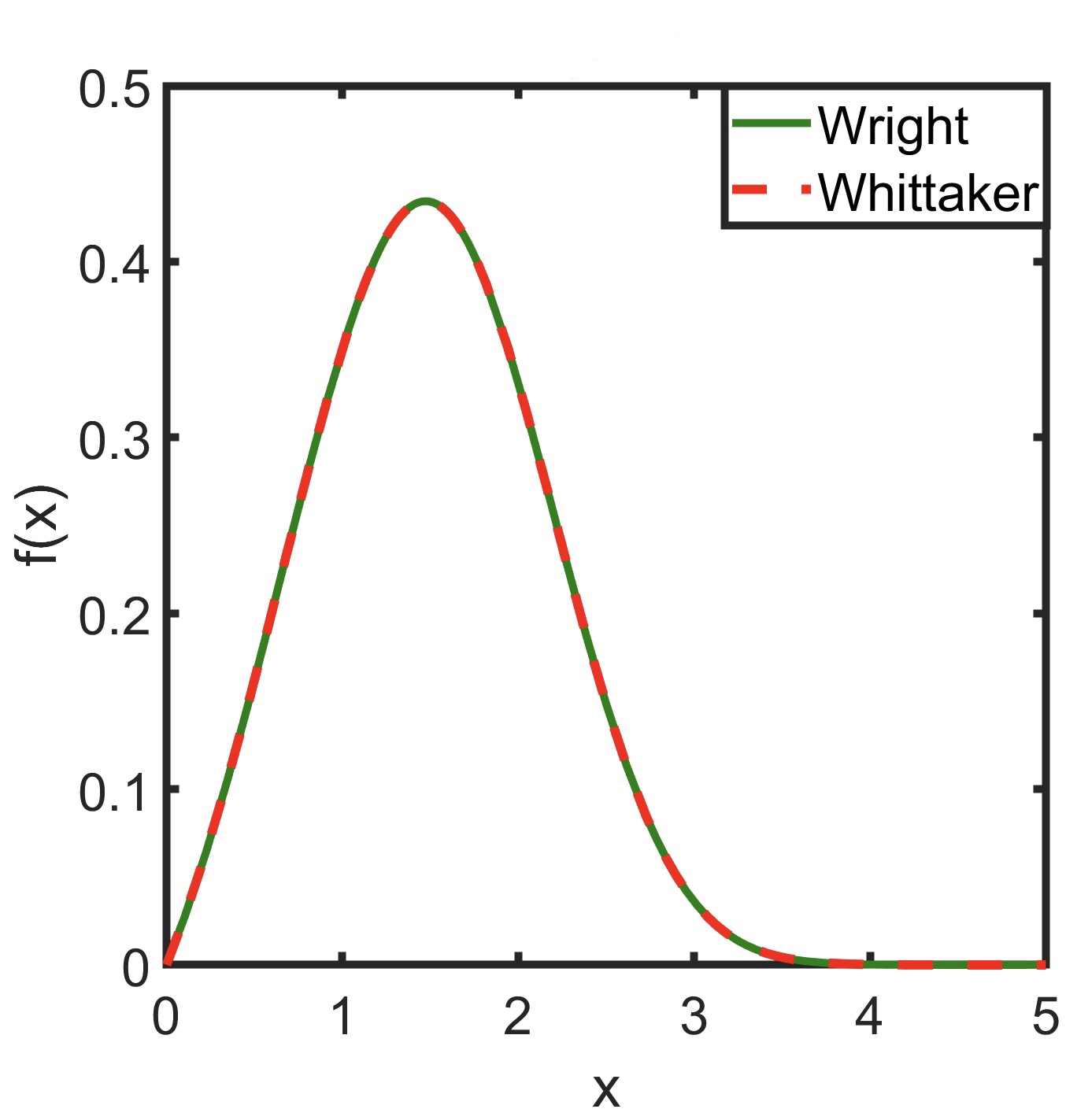}
%\includegraphics[width=.40\textwidth,angle=+0.0]{whittaker_EqC6.jpg}
%\includegraphics[width=.38\textwidth,angle=+0.0]{whittaker_EqC7.jpg}
%\vspace{-0.2truecm}
	\caption{Wright and Whittaker representations
	for (\ref{C6}) (left) and for (\ref{C7}) (right).}
	\label{fig:C6-C7}
\end{figure} %%%%%%%%%%%%%%%%%%%%%

Indeed, replacing $x$ by $x^{-2/3}$ in (\ref{C7}),
we obtain the correct result for the Stankovi\'{c} representation
checked as usual for the identity between the Wright and Whittaker
 representations, see Fig.~\ref{fig:4}, %%% 4, %%%%
\begin{equation}
\label{eq:7.2.12}
W_{-\frac{2}{3}, 0} (- x^{-\frac{2}{3}}) =  \sqrt{\frac{3}{\pi}}
\exp\left(- \frac{2}{27x^2}\right)  {\mathcal W}_{\frac{1}{2}, \frac{1}{6}}
\left( \frac{4}{27x^2}\right).
 \end{equation}

%%%%%%%%%%%%%% Figure 4 %%%%%%%%% 
\begin{figure}[h]
	\centering
	\includegraphics[width=.40\textwidth,angle=+0.0]{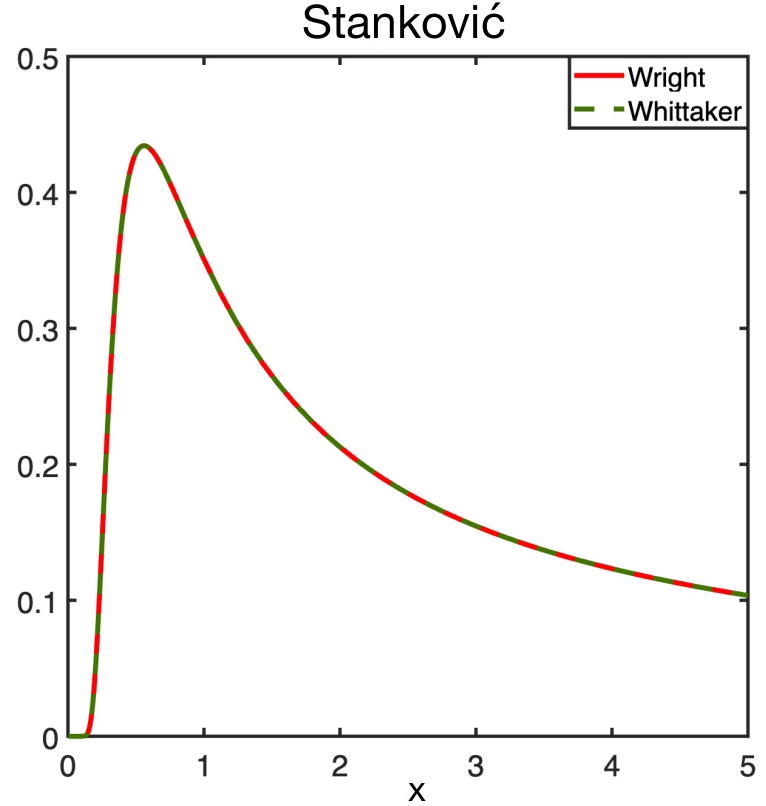}
		%\includegraphics[width=.40\textwidth,angle=+0.0]{eq1dot12.png}
	%% 	\vskip -0.2truecm %%%
		\caption{The corrected Stankovi\'{c} identity
		between the Wright and Whittaker functions for $0\le x\le 5$.}
\label{fig:4} %%%%
 	\end{figure} %%%%%%%%%%%%%%%%%%

%%%%%%%%% Section 3 %%%%%%%%%%%%%%%%%%%%%

\section{A table of special evaluations of the Wright function}\label{sec:3}

\setcounter{section}{3} \setcounter{equation}{0}

In this section we present a list of special evaluations of the Wright function $W_{-\nu,\mu}(\pm x)$ for certain rational values of $\nu$ satisfying $0<\nu<1$ and $0\leq\mu\leq 1$, where
\begin{equation} 
W_{-\nu,\mu}(\pm x)=\sum_{n=0}^\infty \frac{(\pm x)^n}{n! \g(-\nu n+\mu)}\,,
\label{Wright-Paris}
\end{equation}   
with $x>0$.
The functions ${\cal M}_{\kappa,\mu}(x)$ and ${\cal W}_{\kappa,\mu}(x)$ denote the Whittaker functions, $J_\nu(x)$ and $K_\nu(x)$ are the usual Bessel functions and  $_pF_q(x)$ is the generalised hypergeometric function.

The method of proof of these results is the same in each case. An example of the proof when $\nu=2/3$ is supplied in Appendix A.

\bigskip %%% 

{\bf 3.1.} %%% (2.1} %%
\ \ The case $\nu=1/2$,\ $X=x^2/4$:
\begin{eqnarray*}
&&W_{-1/2,\,0}(\pm x)=\mp \frac{X^{1/2}e^{-X}}{\sqrt{\pi}},\\
&&W_{-1/2,1/4}(+x)=\frac{1}{\sqrt{\pi}} X^{-1/4}e^{-X/2}\bl\{{\cal W}_{1/2.1/4}(X)-\frac{\sqrt{\pi}}{\g(\f{3}{4})} 
{\cal M}_{1/2,1/4}(X)\br\},\\
&&W_{-1/2,1/4}(-x)=\frac{1}{\sqrt{\pi}} X^{-1/4}e^{-X/2}{\cal W}_{1/2.1/4}(X),\\
&&W_{-1/2,1/2}(\pm x)=\frac{e^{-X}}{\sqrt{\pi}},\\
&&W_{-1/2,3/4}(+x)=\frac{1}{\sqrt{\pi}} X^{-1/4}e^{-X/2}\bl\{{\cal W}_{0.1/4}(X)+\frac{\sqrt{\pi}}{\g(\f{5}{4})} 
 {\cal M}_{0,1/4}(X)\br\},\\
&&W_{-1/2,3/4}(-x)=\frac{1}{\sqrt{\pi}} X^{-1/4}e^{-X/2}{\cal W}_{0.1/4}(X),\\
&&W_{-1/2,1}(\pm x)=\mp\frac{1}{\sqrt{\pi}}X^{-1/4}e^{-X/2} {\cal W}_{-1/4.1/4}(X)
+2 \binom{1}{0}.
\end{eqnarray*}
The last entry can also be expressed more simply as an error function, namely 
$$
W_{-1/2,1}(\pm x)=1\pm \mbox{erf}\,\sqrt{X}.
$$

\vskip 0.2truecm %%%  

{\bf 3.2.} %% 2.2 %%
 \ The case $\nu=1/3$,\ $X=2(x/3)^{3/2}$:
\begin{eqnarray*}
&&W_{-1/3,\,0}(+x)=-X/2 \{J_{-1/3}(X)+J_{1/3}(X)\}=-3^{-1/3}x\mbox{Ai}(-3^{-1/3}x),\\
&&W_{-1/3,\,0}(-x)=\frac{\sqrt{3}}{2\pi} X\,K_{1/3}(X)=3^{-1/3}x\mbox{Ai}(3^{-1/3}x),\\
&&W_{-1/3,1/3}(+x)= (X/2)^{2/3} \{J_{-2/3}(X)-J_{2/3}(X)\}=-3^{1/3}\mbox{Ai}'(-3^{-1/3}x),\\
&&W_{-1/3,1/3}(-x)= \frac{\sqrt{3}}{\pi} (X/2)^{2/3} K_{2/3}(X)=-3^{1/3}\mbox{Ai}'(3^{-1/3}x),\\
&&W_{-1/3,2/3}(+x)= (X/2)^{2/3} \{J_{-2/3}(X)+J_{2/3}(X)\}=3^{2/3} \mbox{Ai}(-3^{-1/3}x),\\
&&W_{-1/3,2/3}(-x)= \frac{\sqrt{3}}{\pi} (X/2)^{1/3} K_{1/3}(X)=3^{2/3} \mbox{Ai}(3^{-1/3}x),\\
&&W_{-1/3,1}(\pm x)=1\pm \frac{x}{\g(\f{2}{3})}\,{}_1F_2(\f{1}{3}; \f{2}{3}, \f{4}{3};\mp X^2/4)+\frac{x^2}{2\g(\f{1}{3})}\,{}_1F_2(\f{2}{3};\f{4}{3},\f{5}{3};\mp X^2/4).
\end{eqnarray*}

\vskip 0.2truecm %%% %\noindent % 

{\bf 3.3.} %% 2.3 %% 
\ \ The case $\nu=2/3$,\ $X=4x^3/27$:
\begin{eqnarray*}
&&W_{-2/3,\,0}(+x)=-\frac{1}{2\sqrt{3\pi}} e^{X/2} {\cal W}_{-1/2, 1/6}(X),\\
&&W_{-2/3,\,0}(-x)=\sqrt{\frac{3}{\pi}} e^{-X/2} {\cal W}_{1/2, 1/6}(X),\\
&&W_{-2/3,1/3}(+x)= \frac{2^{-4/3}}{\sqrt{3\pi}}\,e^{X/2} X^{-1/3} \,{\cal W}_{-1/2,1/6}(X),\\
&&W_{-2/3,1/3}(-x)=2^{-1/3} \sqrt{\frac{3}{\pi}}\,e^{-X/2} X^{-1/3}\,{\cal W}_{1/2,1/6}(X),\\
&&W_{-2/3,2/3}(+x)= 2^{-2/3}\sqrt{\frac{3}{\pi}} \,e^{X/2} X^{-1/6}\, {\cal W}_{0, 1/3}(X),\\
&&W_{-2/3,2/3}(-x)=2^{-2/3}\sqrt{\frac{3}{\pi}} \,e^{-X/2} X^{-1/6}\, {\cal W}_{0, 1/3}(X),\\
&&W_{-2/3,1}(\pm x)= 1+\frac{2^{-1/3}x}{\sqrt{\pi}}\bl\{\pm \frac{\g(\f{5}{6})}{\g(\f{2}{3})}\,{}_2F_2(\f{1}{3},\f{5}{6};\f{2}{3},\f{4}{3};\pm X)\\
&&\hspace{3cm}-\, X^{1/3} \frac{\g(\f{1}{6})}{4\g(\f{1}{3})}\,{}_2F_2(\f{2}{3},\f{7}{6};\f{4}{3}, \f{5}{3};\pm X)\br\}. \end{eqnarray*}

The cases $\nu=1/4$ and $\nu=3/4$, with $\mu=0, 1/4, 1/2, 3/4, 1$ do not yield any special function representations.
They are found to involve generalised hypergeometric functions of the type
${}_0F_2(-X)$, ${}_1F_3(-X)$ and ${}_2F_3(-X)$, where $X=(x/4)^4$, and so are not included here.

%%%%%%%%%%%%%%%% Section 4 %%%%%%%%%%%%%%%%%%%

\section{The Laplace transform pair occurring in time fractional diffusion processes and the four sisters} \label{sec:4}

\setcounter{section}{4} \setcounter{equation}{0}

In the fractional processes defined by  the partial differential equation of type (\ref{Fractional PDE})
the Wright function of the second kind {\bf is involved}  in the following Laplace transform pair
with $x, t >0$ and $0<\nu<1, \, \mu\ge 0$:
\bee
\label{eq:LT-Stankovic}
\begin{array}{ll}
L^{-1}\left[ s^{-\mu}\, \e^{-x s^\nu}
  \right]
&:= \frac{1}{2\pi i}\int_{Br}\!\! \dfrac{\e^{st-x s^\nu}}{s^\mu}\, ds
= {\ds  t^{\mu-1}} {\ds \sum_{n=0}^\infty}
\dfrac{(-x/t^\nu)^n}{n!\Gamma(- \nu n+\mu)} \\
&:= t^{\mu-1}\, W_{-\nu. \mu}(-x/t^\nu)\,,
\end{array}
\ee
where {\bf $\Re \{ s \}>0$} and Br  stands for the Bromwich path in the complex plane, namely an
infinite line parallel to the imaginary axis cutting the positive real axis
to the right of the branch cut on the negative real axis.
For the reader's convenience, we prove the related inversion of the Laplace transform pair
(\ref{eq:LT-Stankovic}) in Appendix B, where both heuristic and rigorous demonstrations are given.
We note that this pair can also be derived from the 1970 paper by Stankovi\'{c} \cite{Stankovic 70}. 
Other Laplace transform pairs related to $s^{-\mu}\, \exp({-x s^\nu})$ with $x=1$
can be found in the recent article by Apelblat and Mainardi \cite{Apelblat-Mainardi SI21}.

Mainardi and Consiglio  \cite{Mainardi-Consiglio SI20} have utilized the  Laplace transform pair
(\ref{eq:LT-Stankovic}) to define the so-called `four sisters'
relevant in time-fractional diffusion and related to the Mainardi auxiliary functions.
These four functions are obtained for any $\nu \in (0,1)$,
with $\mu=0$, $\mu=1-\nu$, $\mu=\nu$ and $\mu=1$, and
are the natural generalization of  {\it the three sisters functions} obtained for $\nu=1/2$
that henceforth we recall for the reader's convenience  (for more detail, see Appendix A of
 \cite{Mainardi-Consiglio SI20}).
The character of the sisters, because of their inter-relations, was put forward by one of us (F.M.)
in his lecture notes on Mathematical Physics \cite{Mainardi DIFFUSION2019}.
The three sisters, being  related to the fundamental solutions of the standard diffusion equation
obtained from (\ref{Fractional PDE}) with $\beta=1$,  read with their Laplace transforms
\begin{equation} 
     \begin{array} {ll}
   \phi (x,t) &= \textrm{erfc}\left(\dfrac{x}{2\sqrt{t}}\right)\,,
\\ \\
 \psi(x,t) &=\dfrac{x}{2\sqrt{\pi}} \, t^{-3/2} \, \textrm{e}^{-x^{2}/4t}\,,
 \\ \\
 \chi(x,t)&=\dfrac{1}{\sqrt{\pi}} \, t^{-1/2} \, \textrm{e}^{-x^{2}/4t} \,,
\end{array}
\end{equation}
   \vskip -0.3truecm %%
\begin{equation} 
 L[{\phi}(x,t)] = \dfrac{\e^{-x \sqrt{s}}}{s}\,,
\quad
 L[{\psi}(x,t) ]= \e^{-x \sqrt{s}}\,,
 \quad
 L[{\chi}(x,t)] = \dfrac{\e^{-x \sqrt{s}}}{\sqrt{s}}\,.
\end{equation} 
Then, based on the Laplace transform pair (\ref{eq:LT-Stankovic}) we can get for
$\nu=1/2$ and $\mu =1, 0, 1/2$ the representations and plots
of the three sisters for $x,t >0$ in terms of the Wright functions:
\bee
\begin{array}{ll}
\phi(x,t)  &= W_{-1/2, 1} (-x/t^{1/2}),
\\  \\
\psi(x,t) &=  t^{-1} W_{-1/2, 0} (-x/t^{1/2}),
\\  \\
\chi(x,t) &= t^{-1/2} W_{-1/2, 1/2} (-x/t^{1/2}).
\end{array}
\ee

%%%%%%%%%%% Figure 5 %%%%%%%%%%%%%%%%%
\begin{figure}[h!]
	\centering
	\includegraphics[width=.35\textwidth,angle=+0.0]{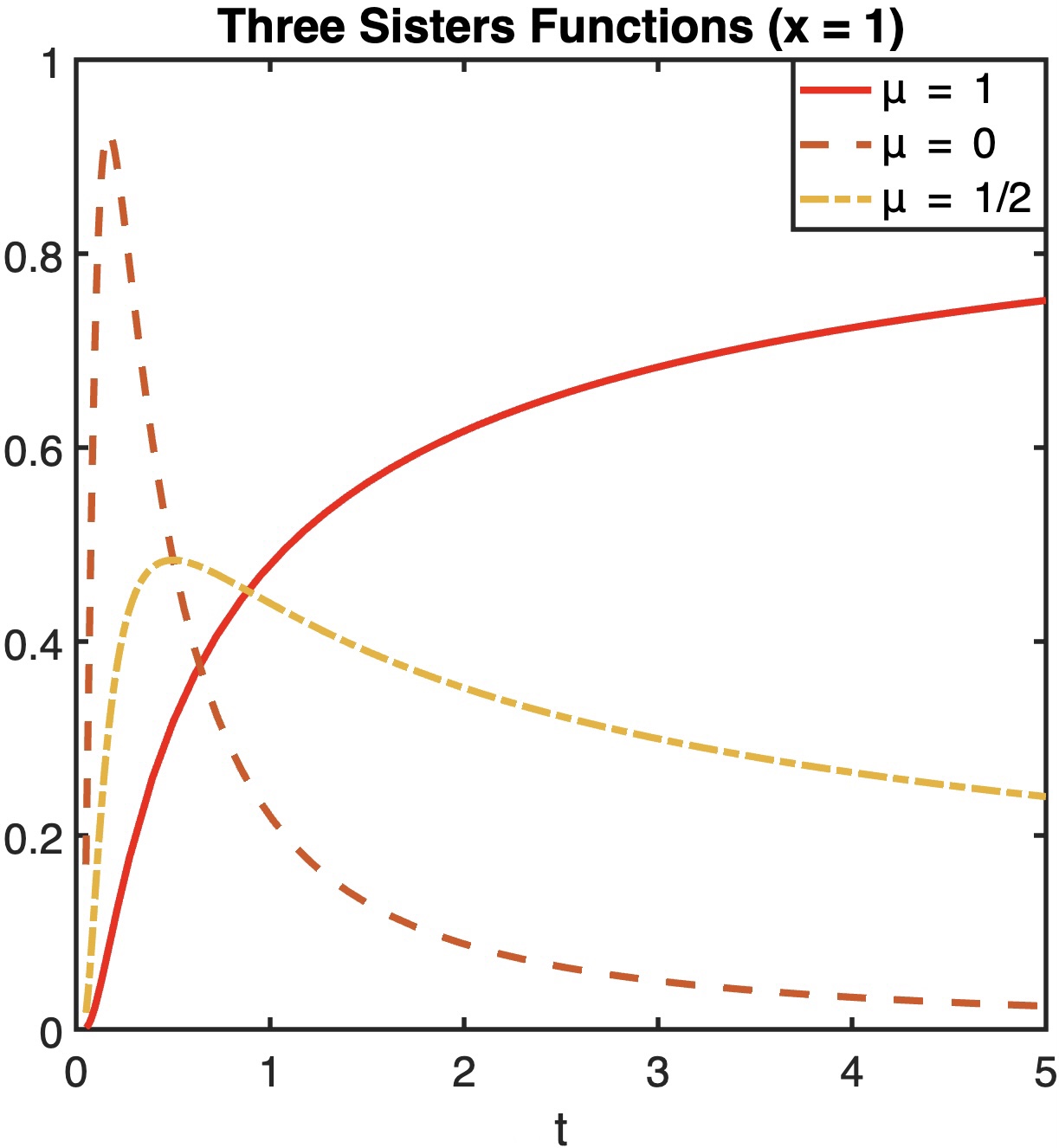}
\hskip 0.5truecm
	\includegraphics[width=.35\textwidth,angle=+0.0]{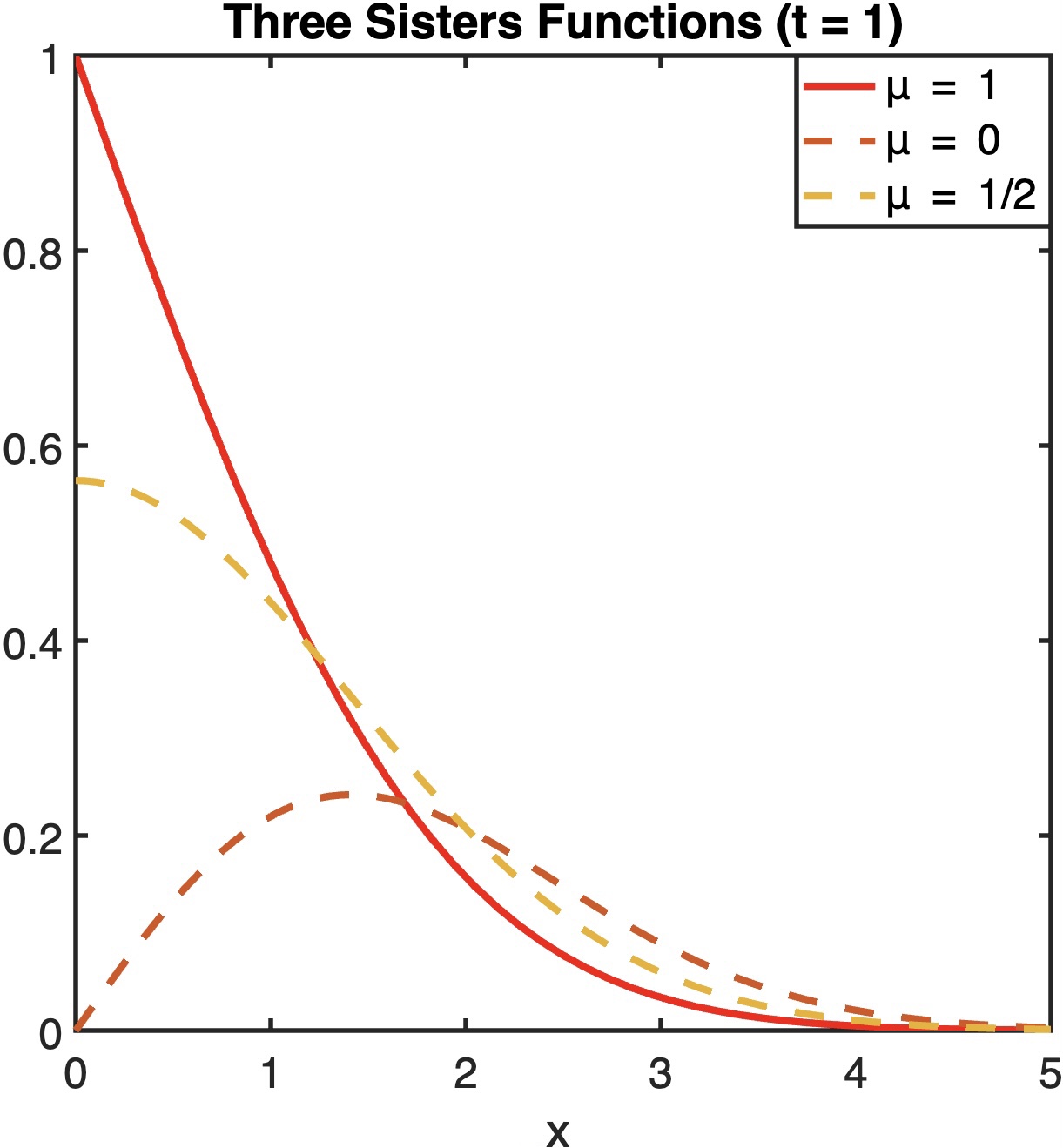}
\caption{The three sisters  versus $t$ at $x=1$ (left) and versus $x$ at $t=1$ (right).} 
	\label{fig:threesisters}
\end{figure} %%%%%%%%%%%%%%%%%%%%%%%%%
%%%
We easily recognize from Section 3.1 the representations of the three sisters  in terms of the Whittaker functions
 by putting $X=x^2/(4t)$ and multiplying by $t^{-\mu}$ accordingly.

In \cite{Mainardi-Consiglio SI20}, the four sisters  were obtained  for  $\nu=1/4, 1/2, 3/4$
and  were plotted versus $x$ at fixed time ($t=1$) and versus $t$ at fixed space ($x=1$).
Then, in Figs. \ref{fig:foursisters_nu1over3}, \ref{fig:foursisters_nu2over3},
we present the new plots  of the four sisters corresponding to  $\nu=1/3$ and $\nu=2/3$.
We can get their representations in terms of the Whittaker functions by using the results of Section 2.2 and 2.3, respectively.

%%%%%%%%%%%% Figure 6 %%%%%%%%%%%%%%%%%%%%%%%%
\begin{figure}[h]
	\centering
	\includegraphics[width=.40\textwidth,angle=+0.0]{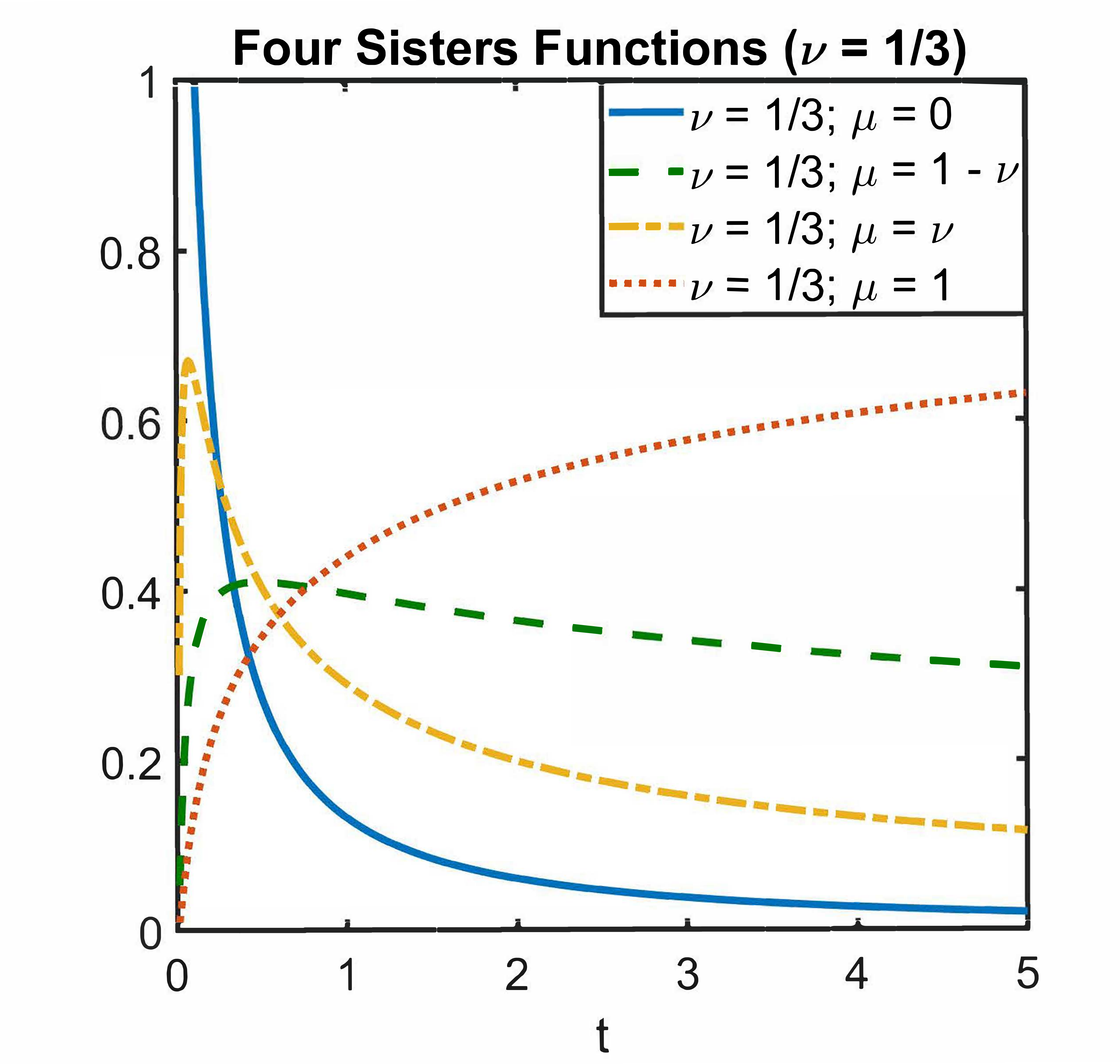}
    \includegraphics[width=.40\textwidth,angle=+0.0]{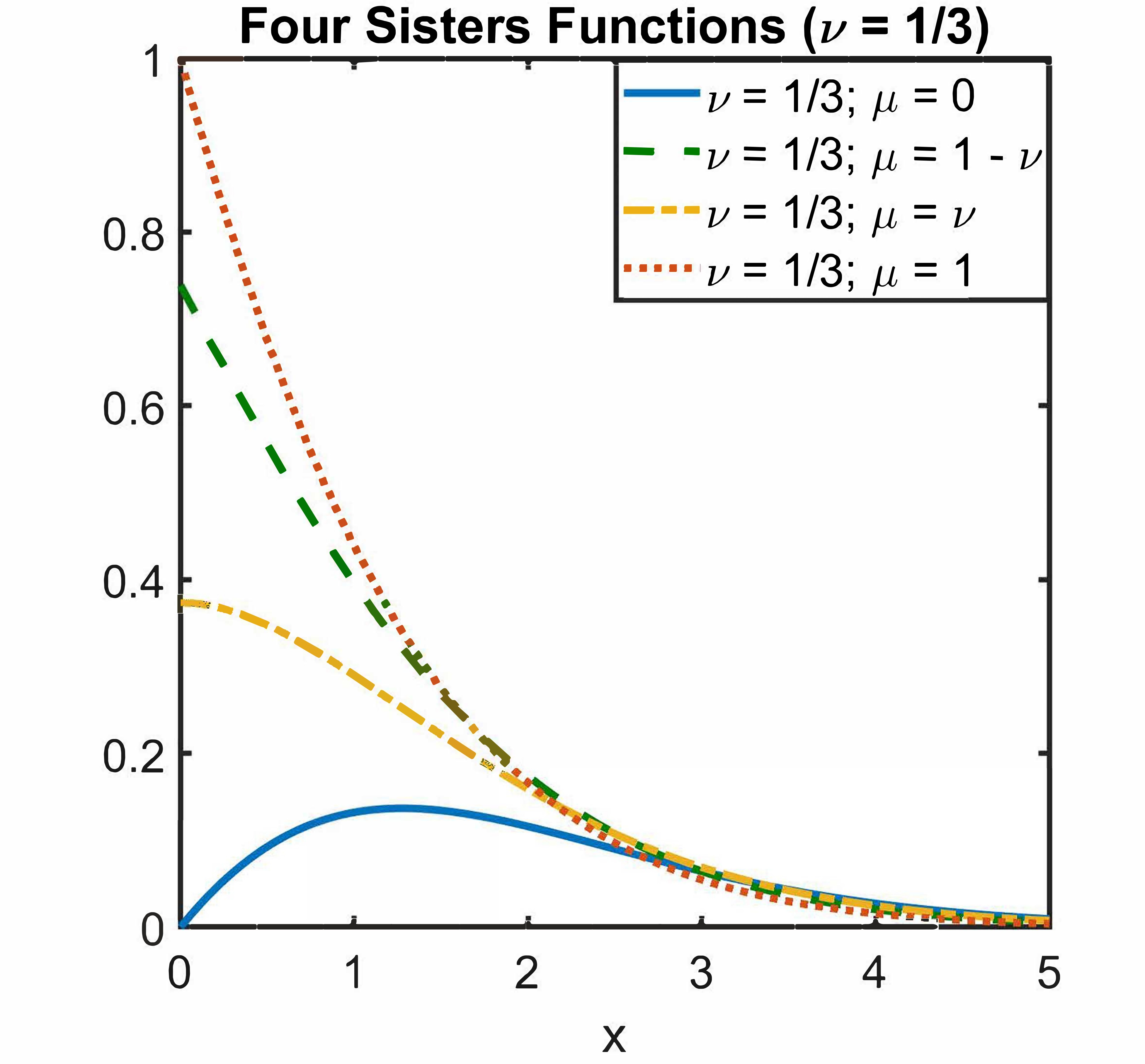}
	\caption{The four sisters for $\nu=1/3$, versus $t$ at $x=1$ (left) and versus $x$ at t=1 (right).}
	\label{fig:foursisters_nu1over3}
\end{figure} %%%%%%%%%%%%%%%%%%%%%

%%%%%%%%%%%%%% Figure 7 %%%%%%%%%%%%%%%%%%%%%
\begin{figure}[h]
	\centering
	\includegraphics[width=.40\textwidth,angle=+0.0]{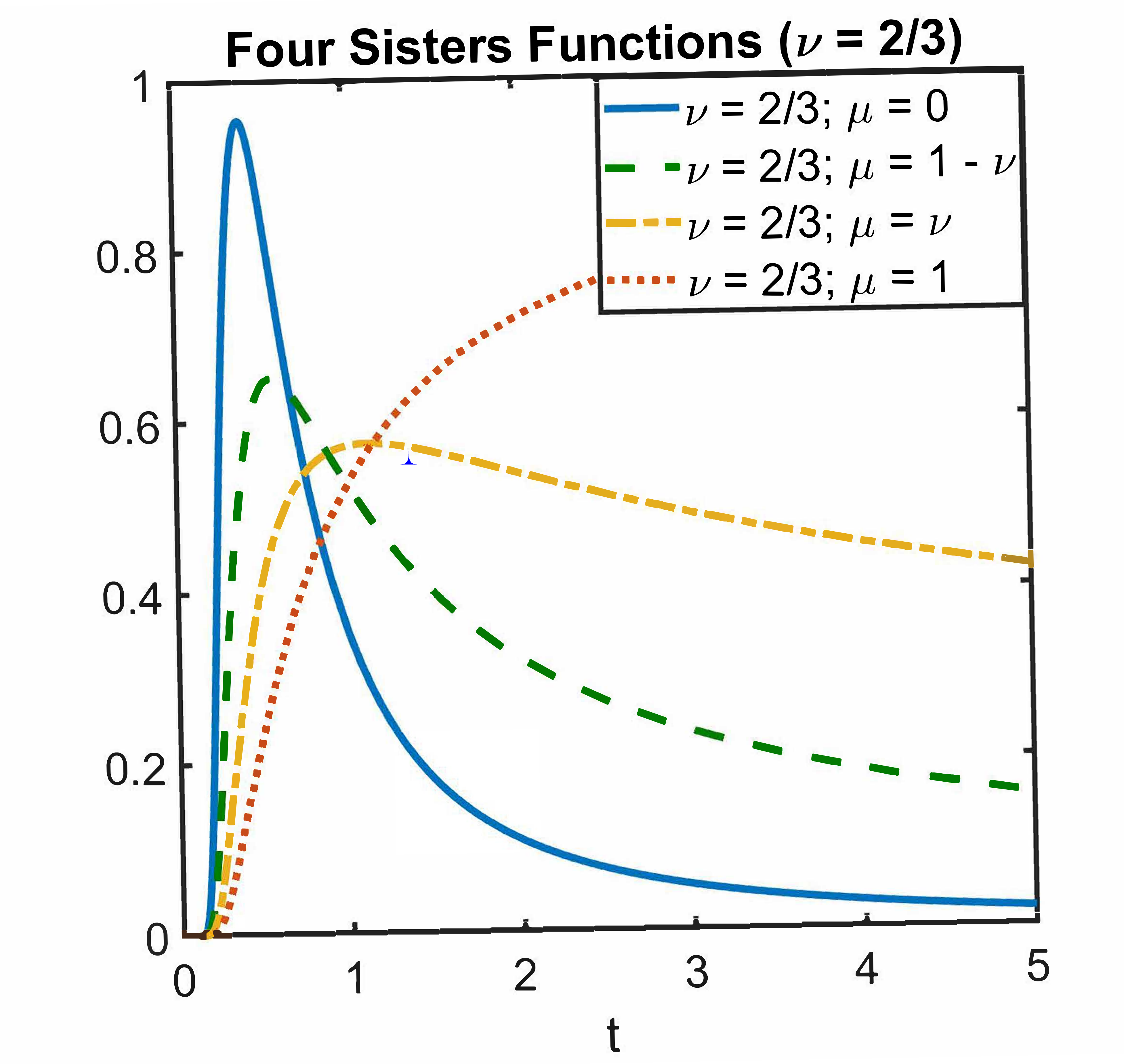}
    \includegraphics[width=.40\textwidth,angle=+0.0]{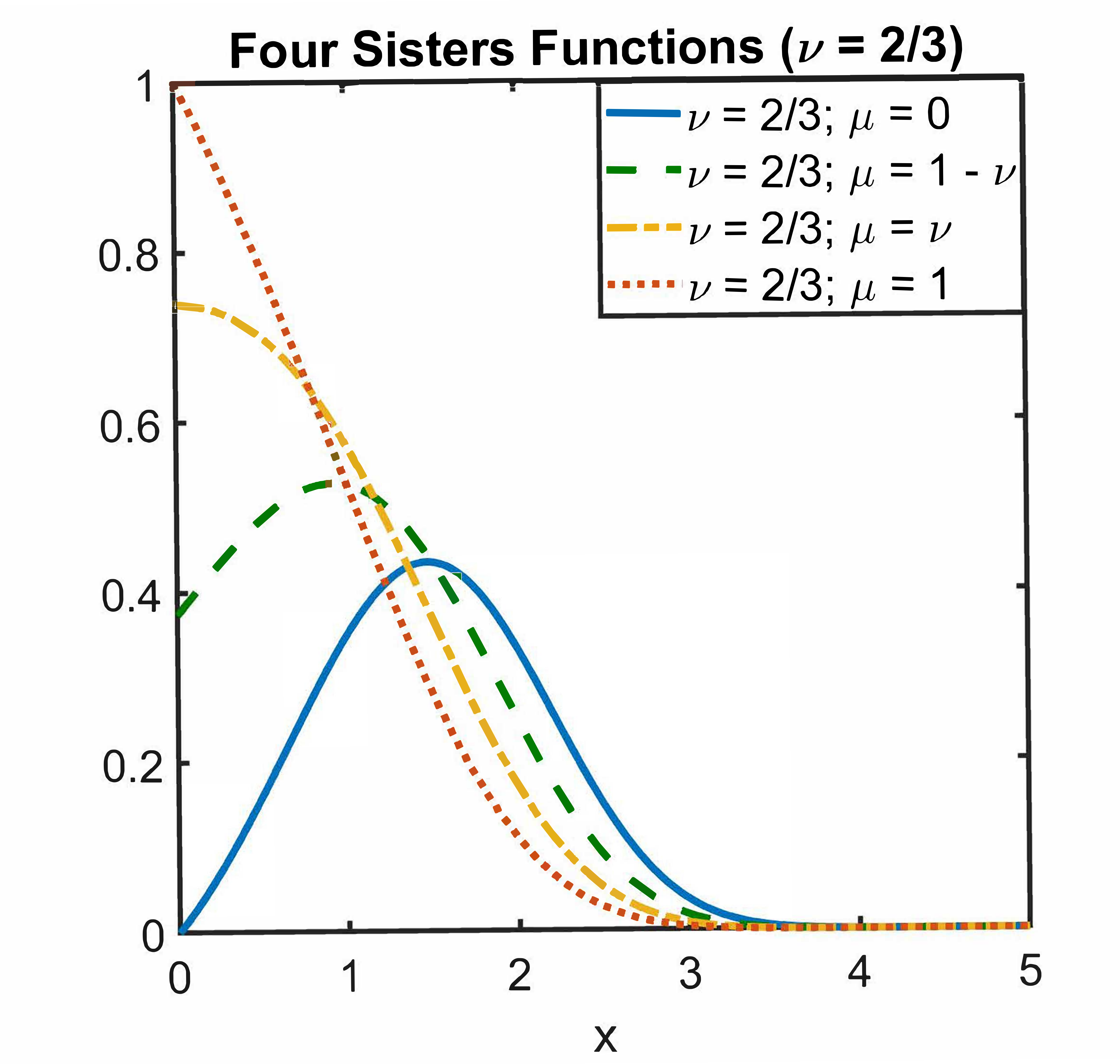}
	\caption{The four sisters for $\nu=2/3$,
 versus $t$ at $x=1$ (left) and versus $x$ at $t=1$  (right).}
	\label{fig:foursisters_nu2over3}
\end{figure} %%%%%%%%%%%%%%%%%%%%%%%%%%%%%%

One sister, namely that corresponding to the left-hand display of Fig.
\ref{fig:foursisters_nu2over3}   with $\nu=2/3$ and $\mu=1-\nu=1/3$,
concerns the case considered by Humbert in 1945  \cite{Humbert 45},
which according to (\ref{eq:LT-Stankovic}) is given by the Laplace transform pair
\bee
\label{LT-Humbert-Wright}
L^{-1}\left[s^{-1/3}\, \exp({- s^{2/3}})\right] =t^{-2/3}\, W_{-{2/3}, {1/3}}(-1/t^{2/3}) \,.
\ee
However, Humbert  %%%%%%%%%%%%%%%%%%%%
\footnote{Humbert used the so-called Carson transform that is the
standard Laplace transform multiplied for $p$, where $p$ stands for  our $s$.}, %%%%%%%
ignoring the Wright function, provided without proof the inverse of the Laplace transform in terms of the Whittaker function as follows, see p.124 in \cite{Humbert 45},
$$ 
L^{-1}\left[s^{-1/3}\, \exp({- s^{2/3}})\right]
 = -\dfrac{1}{4}\, \sqrt{\dfrac{3}{\pi}}\,  \e^{-2/(27t^2)}\,
{\mathcal W}_{-1/2,-1/6}\left(- \dfrac{4}{27t^2}\right),
$$
which is surely wrong.

Indeed, from Section 3.3 we have with
$$ 
x \to 1/t^{2/3}, \quad \, X= 4x^3/27\to 4/(27t^2), 
$$
the expression
\bee 
\label{LT-Humbert-Whittaker}
 t^{-2/3}\, W_{-{2/3}, {1/3}}(-1/t^{2/3})
 = \dfrac{3}{2} \,  \sqrt{\dfrac{3}{\pi}}\, \e^{-2/(27t^2)}\,
 {\mathcal W_{1/2,1/6}}\left( \dfrac{4}{27t^2}\right).
 \ee
The plot in Fig.~\ref{fig:8} % 8 % 
 for $0<t\le 5$ demonstrates the correctness of the expression
(\ref{LT-Humbert-Whittaker}) in terms of Wright and Whittaker functions.
Furthermore,  this plot is seen to be equivalent to that of the sister
for $\nu=2/3$ and $\mu = 1-\nu = 1/3$ in  the left-hand display of 
Fig. \ref{fig:foursisters_nu2over3}, as expected.

%%%%%%%%%%%%%%% Figure 8 %%%%%%%%%%%%%%%%%
\begin{figure}[h!]
	\centering
		\includegraphics[width=.40\textwidth,angle=+0.0]{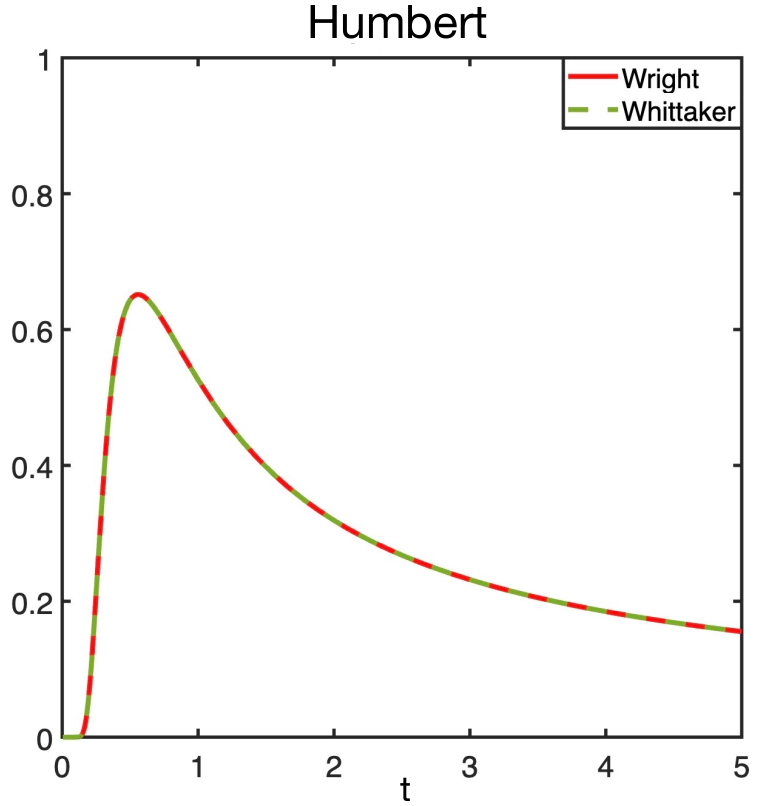}
		\vspace{-0.1truecm} %%
	\caption{The correct  result from the Laplace transform by Humbert for $0< t\le 5$}.
\label{fig:8}%%%
	\end{figure} %%%%%%%%%%%%%%%%%%%%%%%%%

%%%%%%%%%%%%%%%%%% Section 5 %%%%%%%%%%%%%%%%%%%%%%%%%

 \section{Conclusions} \label{sec:5}
 
 A list of evaluations of the Wright function ${\cal W}_{-\nu,\mu}(\pm x)$ of real argument has been given for certain rational values of $\nu$ satisfying $0<\nu<1$ and $0\leq\mu\leq 1$ in terms of familiar special functions, such as the Whittaker, Airy and Bessel functions. 
 Two erroneous representations existing in the literature have been pointed out and corrected.

 The Laplace transform pair occurring in time fractional diffusion processes is established and is employed to define the so-called ``four sisters''.
Plots of these latter functions and other examples of the Wright function are presented to give a graphical indication of their behaviour.

All the plots presented in the paper have been realised using MATLAB. When evaluating power series, it is necessary to employ a sufficiently large number of terms, checking that with fewer terms the results lie within the margin of error.
The Whittaker functions have been plotted using the routines provided by
MATLAB itself whereas the Wright functions have been plotted using their series representation.
 It should be noted that only recently has it been recognised to be advantageous to produce efficient numerical methods for the Wright functions of the second kind, as outlined by Aceto and Durastante \cite{Aceto_2022}.

%%%%%%%  APPENDICES  %%%%%%%%%%%%%%%%%%%%%%%%%%%%%%%%%%%%%%%%%%% 

\section*{Appendix A: Proof of the case $\nu=2/3$ in \S 3.3}

\setcounter{section}{1} \setcounter{equation}{0}
\renewcommand{\theequation}{\Alph{section}.\arabic{equation}} %%%%%%%%%%%%%%

In this appendix we prove the Whittaker function evaluations given in Section 3.3.
The proofs of the other results in Section \ref{sec:3} 
follow the same procedure and are consequently omitted.
We consider the case $\nu=2/3$ for $0\leq \mu\leq 1$ and write $W_{-2/3,\mu}(\pm x)$ in the form
$$
W_{-2/3,\mu}(\pm x) =\frac{1}{\pi}\sum_{n=0}^\infty \frac{(\pm x)^n}{n!}\,\g(\f{2}{3}n+1-\mu) \sin \pi(\f{2}{3}n+1-\mu).
$$
Replacement of $n$ by $3m+j$, $j=0, 1, 2$ then produces
\begin{equation*}
\begin{array}{ll}
W_{-2/3,\mu}(\pm x)&=
{\ds \frac{\sin \pi(1-\mu)}{\pi}\sum_{m=0}^\infty \frac{(\pm x)^{3m}}{(3m)!} \g(2m\!+\!1\!-\!\mu)}\\
&\pm {\ds \frac{x \sin \pi(\f{5}{3}-\mu)}{\pi}\sum_{m=0}^\infty \frac{(\pm x)^{3m}}{(3m+1)!} \g(2m\!+\!\f{5}{3}\!-\!\mu)}\\
&+
{\ds \frac{x^2 \sin \pi(\f{7}{3}-\mu)}{\pi}\sum_{m=0}^\infty \frac{(\pm x)^{3m}}{(3m+2)!}\,\g(2m\!+\!\f{7}{3}\!-\!\mu)}.
\end{array}
\end{equation*}

Using the reflection and triplication formulas for the gamma function
\[
\g(2a)=\frac{2^{2a-1}}{\sqrt{\pi}}\,\g(a) \g(a+\fs),\qquad \g(3a)=\frac{3^{3a-\frac{1}{2}}}{2\pi}\,\g(a)\g(a+\f{1}{3})\g(a+\f{2}{3})
\]
and introducing of the new variable $X=4x^3/27$, we find
\begin{eqnarray}
 & & W_{-2/3,\mu}(\pm x)
 = %%% &=& %% to avoid hfull ..% 
 \frac{2^{1-\mu}}{\sqrt{3\pi}} \bl\{\sin \pi(1-\mu) \sum_{m=0}^\infty (\pm X)^m \,\frac{\g(m+\fs-\fs\mu) \g(m+1-\fs\mu)}{\g(m+\f{1}{3}) \g(m+\f{2}{3})\,m!} \nonumber\\
& & \hspace*{0.5cm} %% 
\pm \, 
X^{1/3} \sin \pi(\f{5}{3}-\mu) \sum_{m=0}^\infty (\pm X)^m \frac{\g(m+\f{5}{6}-\fs\mu)\g(m+\f{4}{3}-\fs\mu)}{\g(m+\f{2}{3})\g(m+\f{4}{3})\,m!}\hspace{1cm}\nonumber\\
& & \hspace*{0.5cm} %% 
+ \, X^{2/3} \sin \pi(\f{7}{3}-\mu) \sum_{m=0}^\infty (\pm X)^m \frac{\g(m+\f{7}{6}-\fs\mu) \g(m+\f{5}{3}-\fs\mu)}
{\g(m+\f{4}{3}) \g(m+\f{5}{3})\,m!}\br\}. \nonumber \\ %%%
\label{a1}
\end{eqnarray}

{\it The case $\mu=1/3$\/.}
\vspace{0.2cm} %%%

When $\mu=1/3$, we have from (\ref{a1})
\[
W_{-2/3,1/2}(\pm x)=\frac{2^{-1/3}}{\sqrt{\pi}}\bl\{\sum_{m=0}^\infty \frac{(\pm X)^m \g(m+\f{5}{6})}{\g(m+\f{2}{3})\,m!} \mp X^{1/3} \sum_{m=0}^\infty \frac{(\pm X)^m \g(m+\f{7}{6})}{\g(m+\f{4}{3})\,m!}\br\}
\]
\[
=\frac{2^{-1/3}}{\sqrt{\pi}}\bl\{\frac{\g(\f{5}{6})}{\g(\f{2}{3})}\,{}_1F_1(\f{5}{6};\f{2}{3};\pm X)\mp X^{1/3} \frac{\g(\f{7}{6})}{\g(\f{4}{3})}\,{}_1F_1(\f{7}{6};\f{4}{3};\pm X)\br\}
\hspace*{0.4cm} %% 
\]
\begin{equation}\label{a2}
\ \ = \frac{2^{-4/3}}{\sqrt{3\pi}}\bl \{\frac{\g(\f{1}{3})}{\g(\f{7}{6})}\,{}_1F_1(\f{5}{6};\f{2}{3};\pm X)\pm X^{1/3} \frac{\g(-\f{1}{3})}{\g(\f{5}{6})}\,{}_1F_1(\f{7}{6};\f{4}{3};\pm X)\br\},
\end{equation}
where ${}_1F_1(a;b;z)$ is the confluent hypergeometric function.
From (\ref {e18a}) we have the following combination of confluent hypergeometric functions expressed in terms 
of the Whittaker function ${\cal W}_{\kappa,\mu}(z)$:
\[
\frac{\g(1-b)}{\g(a-b+1)}\,{}_1F_1(a;b;z)+z^{1-b} \frac{\g(b-1)}{\g(a)}\,{}_1F_1(a-b+1;2-b;z)
\]
   \renewcommand{\theequation}{\Alph{section}.\arabic{equation}} %%%%%%%%%%%%%% 
\begin{equation} \label{a3} %%%%
\hspace{3cm} %% 
=e^{z/2} z^{-b/2} {\cal W}_{a/2-a,b/2-1/2}(z)\quad (z\in {\bf C}).
\end{equation}
Inserting the values $a=5/6$ and $b=2/3$ into (\ref{a3}) we find upon taking the upper signs in (\ref{a2}) that
\begin{equation} \label{a4}
W_{-2/3,1/3}(x)=\frac{2^{-4/3}}{\sqrt{3\pi}}\,e^{X/2} X^{-1/3} {\cal W}_{-1/2,1/6}(X),
\end{equation}
where we have made use of the fact that ${\cal W}_{\kappa,-\mu}(z)={\cal W}_{\kappa,\mu}(z)$.

Taking the lower sign in (\ref{a2}), we have upon application of the well-known Kummer transformation 
${}_1F_1(a;b;z)=e^z {}_1F_1(b-a;b;-z)$
\[  %% splitted ! %%%
W_{-2/3,1/3}(-x)=2^{-1/3} \sqrt{\frac{3}{\pi}}\,e^{-X} \bl\{\frac{\g(\f{1}{3})}{\g(\f{1}{6})}
\,{}_1F_1(-\frac{1}{6};
\frac{2}{3};X)
\] %%%%%
\vspace*{-0.2cm} %%%
\[ \hspace*{1cm} %%%
+\, X^{1/3}\,\frac{\g(-\f{1}{3})}{\g(-\frac{1}{6})}\,{}_1F_1(\f{1}{6};\f{4}{3};X)\br\}.
\]
Insertion of the values $a=-1/6$ and $b=2/3$ into (\ref{a3}) then yields
\begin{equation}\label{a5}
W_{-2/3,1/3}(-x)=2^{-1/3} \sqrt{\frac{3}{\pi}}\,e^{-X/2} X^{-1/3} {\cal W}_{1/2,1/6}(X).
\end{equation}

{\it The case $\mu=2/3$\/.}
\vspace{0.2cm} %%%

%% \noindent
Proceeding in the same manner when $\mu=2/3$, we have from (\ref{a1})
\[W_{-2/3,2/3}(\pm x)=\frac{2^{-2/3}}{\sqrt{\pi}}\bl\{\sum_{m=0}^\infty \frac{(\pm X)^m \g(m+\f{1}{6})}{\g(m+\f{1}{3})\,m!}-X^{2/3}\sum_{m=0}^\infty \frac{(\pm X)^m \g(m+\f{5}{6})}{\g(m+\f{5}{3})\,m!}\br\}\]
\[=2^{-2/3}\sqrt{\frac{3}{\pi}} \bl\{\frac{\g(\f{2}{3})}{\g(\f{5}{6})}\,{}_1F_1(\f{1}{6};\f{1}{3};\pm X)
+X^{2/3}\,\frac{\g(-\f{2}{3})}{\g(\f{1}{6})}\,{}_1F_1(\f{5}{6}; \f{5}{3};\pm X)\br\}.\]
With the values $a=1/6$, $b=1/3$ in (\ref{a3}) we find upon taking the upper signs that
\begin{equation}\label{a6}
W_{-2/3,2/3}(x)=2^{-2/3}\sqrt{\frac{3}{\pi}}\,e^{X/2} X^{-1/6}\,{\cal W}_{0,1/3}(X).
\end{equation}
Application of Kummer's transformation to the expression with the lower signs and the values $a=1/6$, $b=1/3$ in (\ref{a3}) then yields
\begin{equation}\label{a7}
W_{-2/3,2/3}(-x)=2^{-2/3}\sqrt{\frac{3}{\pi}}\,e^{-X/2} X^{-1/6}\,{\cal W}_{0,1/3}(X).
\end{equation}

{\it The case $\mu=0$\/.}
\vspace{0.2cm} %% 

Finally, when $\mu=0$ we have from (\ref{a1})
\[
W_{-2/3,0}(\pm x)=\mp\frac{X^{1/3}}{\sqrt{\pi}} \bl\{\sum_{m=0}^\infty \frac{(\pm X)^m \g(m+\f{5}{6})}{\g(m+\f{2}{3})\,m!}\pm X^{1/3} \sum_{m=0}^\infty \frac{(\pm X)^m \g(m+\f{7}{6})}{\g(m+\f{4}{3})\,m!}\br\}
\]
\vspace*{-0.3cm} %%%%
\[
\hspace{2cm}
=\mp\frac{X^{1/3}}{2\sqrt{3\pi}} \bl\{\frac{\g(\f{1}{3})}{\g(\f{7}{6})}\,{}_1F_1(\f{5}{6};\f{2}{3};
\pm X)\pm X^{1/3} \frac{\g(-\f{1}{3})}{\g(\f{5}{6})}\,{}_1F_1(\f{7}{6};\f{4}{3};\pm X)\br\}.
\]
Taking the upper signs and the values $a=5/6$, $b=2/3$ in (\ref{a3}), we find
\begin{equation}\label{a8}
W_{-2/3,0}(x)=-\frac{1}{2\sqrt{3\pi}}\,e^{X/2}\,{\cal W}_{-1/2,1/6}(X),
\end{equation}
and with the lower signs, after application of Kummer's transformation, 
we have with $a=-1/6$, $b=2/3$ in (\ref{a3})
\[
W_{-2/3,0}(-x)\!=\! \sqrt{\frac{3}{\pi}} X^{1/3} e^{-X} \bl\{\frac{\g(\f{1}{3})}{\g(\f{1}{6})}\,{}_1F_1(-\f{1}{6};\f{2}{3};X)+X^{1/3} \frac{\g(-\f{1}{3})}{\g(-\f{1}{6})}\,{}_1F_1(\f{1}{6};\f{4}{3};X)\br\}
\]
\begin{equation}\label{a9}
=\sqrt{\frac{3}{\pi}}\,e^{-X/2}\,{\cal W}_{1/2,1/6}(X).
\end{equation}

The case $\mu=1$ in ({\ref{a1}) results in no cancellations of the gamma functions in the second and third series (the first series vanishing identically). As a consequence, these series are expressible in terms of ${}_2F_2(\pm X)$ functions as stated in Section 2.3, with the result that they will not reduce to Bessel or Whittaker functions.

%%%%%%%%%%%%%%%%%%%%%%%%

\section*{Appendix B: On inversion of a relevant  Laplace transform}

\setcounter{section}{2} \setcounter{equation}{0}
\renewcommand{\theequation}{\Alph{section}.\arabic{equation}} %%%%

In this appendix we evaluate the inverse Laplace transform (LT) (\ref{eq:LT-Stankovic})
whose expression is repeated here for the reader's convenience,
\bee  \label{B.1}
L^{-1}\left[s^{-\mu}\, \exp({-x s^\nu})\right] = t^{\mu-1}\, W_{-\nu. \mu}(-x/t^\nu)\,,
\  \; 0<\nu<1, \ \mu>0\,.
\ee

We first provide a formal justification of this result starting from the following  LT pair reported by Doetsch
for the (extended) Laplace transforms of
\bee \label{B.2}
L\left[Pf \dfrac{t^{-\alpha-1}}{\Gamma(-\alpha)}\right] = s^\alpha,
\ \; \alpha>0, \ \alpha \ne 0, 1, 2, \dots,
\ee
where $Pf$ denotes the Hadamard finite part, see \cite[pp.~63, 318]{Doetsch BOOK74}.
For $\alpha = n$, with $n= 0, 1, 2$, we obtain the Laplace transforms of the generalised functions
$\delta^{(n)}(t)$, that is of the derivatives of the Dirac`delta function'.
Indeed it is known that the standard Laplace transforms of locally integrable functions are required 
to  vanish at infinity. We note that (B.2) can also be derived from the standard LT pair
\bee\label{B.3}
L\left[t^\alpha \right] = \dfrac{\Gamma(\alpha+1)}{s^{\alpha+1}}, \quad  \alpha>-1\,.
\ee

Now we proceed formally expanding the exponential series  term by term as follows taking into account 
the LT pair (B.2) for $t>0$:
\bee \label{B.4}
\begin{array}{ll}
&L^{-1}\left[ {\ds s^{-\mu}\, \exp({-x s^\nu})}\right] =
{\ds \sum_{n=0}^\infty} \dfrac{(-x)^n}{n!}\, L^{-1}\left[ {\ds s^{\nu n -\mu}}\right]\\  %% \\ %
& = {\ds \sum_{n=0}^\infty} \dfrac{(-x)^n}{n!}\,  \,
\dfrac{t^{-\nu n +\mu-1}}{\Gamma(-\nu n+\mu)}
:= {\ds t^{\mu-1}}\, W_{-\nu, \mu}\left(-x/t^\nu\right)\,.
\end{array}
\ee
We note that this result is considered valid for $t>0$; that is excluding the point $t=0$
because when the exponent of the Laplace parameter turns out to be a non-negative integer
the contributions of the corresponding generalised functions vanish.

We now provide a rigorous justification of the result (B.1) using the Bromwich integral formula for the inverse LT followed by use of the Hankel formula for the reciprocal of the gamma function given by
\bee \label{B.5}
\dfrac{1}{\Gamma (z)} 
= \dfrac{1}{2\pi i} \, \int_{Ha}
 \dfrac{\e^{\tau}}{{\tau^z}}\, d\tau \,, \quad  z \in \CC \,,
 \ee 
where $Ha$ denotes the Hankel contour defined as a path that
begins at $\tau =-\infty - ia$ ($a>0$), encircles the  branch
cut that lies along the negative real axis, and terminates at
$\tau= - \infty + ib$ ($b>0$).
The branch cut is necessary when $z$ is non-integer because $\tau^{-z}$ is a multivalued function.
Then, provided $0<\nu<1$, the Bromwich path may be deformed into the Hankel path to yield
\bee\label{B.6}
\begin{array}{ll}
& L^{-1} \left [ s^{-\mu}\, {\ds \e ^{-xs^\nu}} \right]
= \dfrac{1}{2\pi i} \, \int_{Ha}
{\ds \e^{st}}\, \dfrac{\e^{-xs^\nu}}{s^\mu} \, ds \\ \\
&= {\ds \sum_{n=0}^\infty} \dfrac{(-x)^n}{n!}
\dfrac{1}{2\pi i} \, \int_{Ha}\dfrac{ \e^{st}}{s^{\mu-\nu n}} \, ds \\ \\
&= {\ds \sum_{n=0}^\infty} \dfrac{(-x)^n}{n!}\,
t^{\mu -\nu n-1} \dfrac{1}{2\pi i} \, \int_{Ha_-}\dfrac{ \e^{\tau}}{\tau^{\mu-\nu n}} \, d\tau
= t^{\mu-1}\, W_{-\nu, \mu}\left( -x/t^\nu\right)\,.
\end{array}
\ee

%%% moved here !!! %%%%%%%%%%%%%%%%%%%%%%%%%%%%%%%%

\subsection*{Acknowledgements} %%%

The research activity of A.C. and F.M. has been carried out in the framework of the activities of the National Group of Mathematical Physics (GNFM, INdAM).
The activity of A.C., a PhD student at the University of Wuerzburg, is carried out also in the Wuerzburg-Dresden Cluster of Excellence - Complexity and Topology in Quantum Matter (ct.qmat).
\\
We wish to thank the anonymous referees for their helpful suggestions that have improved the presentation of the paper.
% \end{acknowledgements} %%%

 \section*{\small Conflict of interest} %%%%%%%%%%%%%%%%%%%%
 {\small The authors declare that they have no conflict of interest.}
 
 %%%%%%%%%%%%%%%%%%%%%%%%%%%%%%%%%%%%%%%%%%%%%%%%%%%%%%%%%%%

\bigskip  %%%%%%%%%%%%%%%%%%%%%%%%%%%%%%%%%

\end{document}